\documentclass[reqno]{amsart}
\usepackage{hyperref}
\usepackage{color}

\def\N{{\mathbb{N}}}

\def\R{{\mathbb{R}}}

\begin{document}
\title[Multiplier rules]{A generalization of multiplier rules for infinite-dimensional optimization problems}
\author[H. Yilmaz]{Hasan Yilmaz}

\address{Hasan Yilmaz: Laboratoire LPSM UMR 8001 \newline
Universit\'e Paris-Diderot, Sorbonne-Paris-Cit\'e \newline
b\^atiment Sophie Germain, 8 place Aur\'elie Nemours, \newline
75013 Paris, France.}
\email{yilmaz@lpsm.paris}
\date{January, 17, 2019}
\numberwithin{equation}{section}
\newtheorem{theorem}{Theorem}[section]
\newtheorem{lemma}[theorem]{Lemma}
\newtheorem{example}[theorem]{Example}
\newtheorem{remark}[theorem]{Remark}
\newtheorem{definition}[theorem]{Definition}
\newtheorem{corollary}[theorem]{Corollary}
\newtheorem{proposition}[theorem]{Proposition}
\thispagestyle{empty} \setcounter{page}{1}
\begin{abstract}
We provide a generalization of first-order necessary conditions of optimality for infinite-dimensional optimization problems with a finite number of inequality constraints and with a finite number of inequality and equality constraints. Our assumptions on the differentiability of the functions are weaker than those of existing results.    
\end{abstract} 
\maketitle
\vskip3mm
\section{Introduction}
We provide an improvement of first-order necessary conditions of optimality for infinite-dimensional problems under a finite number of inequality constraints and under a finite number of inequality and equality constraints in the form of Fritz John's theorem and Karush-Kuhn-Tucker's theorem. \\
In this paper, we give a proof of multiplier rules by following the same approach as Michel in \cite{PM} p. 510. To prove his result, Michel uses the Brouwer fixed-point theorem that is why, at the first slight, his proof seems specific to finite-dimensional optimisation problems. However, we remark that we can extend this result for infinite-dimensional optimizaton problems. The proof of Michel is explained in detail in \cite{BH}, Appendix B. Another proof of the multiplier rules was established by Halkin in \cite{HH} but his proof is completely different. Indeed, Halkin uses an implicit function theorem with only Fr\'echet differentiable at a point framework instead of the continuously Fr\'echet differentiable framework. The improvement of Michel and Halkin is to replace the assumption of continuously Fr\'echet differentiable on a neighborhood of the optimal solution (see in \cite{PB} Chapter 13 section 2) with the assumptions of the continuity on a neighborhood of the optimal solution and the Fr\'echet differentiability at the optimal solution.   \\ 
Note that there are another way to generalize the assumption of continuous Fr\'echet differentiability by using locally Lipschitzian mappings  e.g. \cite{CL}. The statement of Halkin and Michel is not similar with the statements of locally Lispchitzian. Indeed, in general, a mapping which is Fr\'echet differentiable at a point is not locally Lipschitzian arround this point and conversely a mapping which is locally Lipschitzian arround a point is not Fr\'echet differentiable at this point.\\
In \cite{BL}, Blot gave also a proof of multiplier rules for finite-dimensional optimization problems under only inequality constraints and under inequality and equality constraints. For the problems with inequality constraints, Blot reduced the assumptions of Pourciau in \cite{P} by replacing, at the optimal solution, the Fr\'echet differentiability with the G\^ateaux differentiability. For the problems with inequality and equality constraints, Blot deleted the assumptions of local continuity on a neighborhood for the objective function and for the functions in the constraints of inequality. Therefore, Blot improved the multilplier rules of Michel and Halkin by lightening the assumptions on the continuity of the functions.\\
The main contributions of this paper are as follows.
\begin{itemize}
\item Contrary to Blot and Michel, we do not assume the finiteness of the dimension of the space. Therefore, we extend the main theorems of \cite{BL} and \cite{PM} in infinite-dimensional vector spaces.
\item Moreover, in comparison with Blot's multiplier rules, we replaced the assumptions of Fr\'echet differentiability by the Hadamard differentiability which is weaker in infinite-dimensional vector. For consequently, our assumptions on the differentiability of the functions are weaker than \cite{BL}, \cite{PM} and \cite{HH}.
\end{itemize}
We summarize the content of this paper as follows.\\
In Section 2, we state the main theorems of the paper.\\
In Section 3, we specify the definition of G\^ateaux differentiability and Hadamard differentiability. Besides, we recall a supporting hyperplane theorem and the \\Schauder fixed-point theorem.\\
In Section 4, in order to proof our first-order necessary conditions under inequality constraints, we delete the inactive inequality constraints. Next, we use a supporting hyperplane theorem to find the multipliers.\\
In Section 5, we give a proof of first-order necessary conditions of optimality under inequality and equality constraints. As in Section 4, we delete the inactive inequality constraints. In order to use the supporting hyperplane theorem, we use the Schauder fixed-point theorem. 
\section{Statements of the Main Results} The paper deals with infinite-dimensional optimization problems with a finite list of inequality constraints and with a finite list of inequality and equality constraints.
Let $E$ be a normed vector space, let $\Omega$ be a nonempty open subset of $E$, let $f_i:\Omega \rightarrow \R$ when $i\in\{0,...,m\}$ be functions, let $f:\Omega \rightarrow \R$, $g_i:\Omega \rightarrow \R$ when $i\in\{1,...,p\}$, $h_j: \Omega \rightarrow \R$ when $j\in\{1,...,q\}$ be functions and $m$, $p$ and $q$ are integer number.
We consider the two following problems \[
({\mathcal I})
\left\{
\begin{array}{cl}
{\rm Maximize} & f_0(x)  \\
{\rm subject \;  to} & x \in \Omega \\
\null & \forall i\in\{1,..., m\},\, f_i(x) \ge 0
\end{array}\right.
\] 
and 
\[({\mathcal P})
\left\{
\begin{array}{cl}
{\rm Maximize} & f(x)  \\
{\rm subject \;  to} & x \in \Omega \\
\null & \forall i\in\{1,..., p\},\, g_i(x) \ge 0\\
\null & \forall j\in\{1,..., q\},\, h_j(x) = 0.\\
\end{array}\right.
\]
The main theorems of the paper are the following ones.
\begin{theorem}\label{th21} Let $\hat{x}$ be a solution of $({\mathcal I})$. We assume that the following assumptions are fulfilled.\\
{\rm (i)} For all $i\in \{0,...,m\}$, $f_i$ is G\^ateaux differentiable at $\hat{x}$.\\
{\rm (ii)} For all $i\in \{1,...,m\}$, $f_i$ is lower semicontinuous at $\hat{x}$ when $f_i(\hat{x})>0$.\\
Then, there exist $\lambda_0,..., \lambda_m \in \R_+$ which satisfy the following conditions.\\
\rm{(a)} $(\lambda_0,..., \lambda_m) \ne (0,...,0).$\\
\rm{(b)} For all $i \in \{1,..., m\},\; \lambda_if_i(\hat{x})=0.$\\
\rm{(c)} $\sum_{i=0}^{m} \lambda_iD_Gf_i(\hat{x})=0.\\
$ 
In addition, if we assume that the following assumption is verified  \\
{\rm (iii)} there exists $w\in E$ such that for all $i\in \{1,...,m\}, D_Gf_i(\hat{x})w>0$ when $f_i(\hat{x})=0$\\
then we can take \\
{\rm (d)} $\lambda_0=1$.
\end{theorem}
\begin{theorem}\label{th22} Let $\hat{x}$ be a solution of $({\mathcal P})$. We assume that the following assumptions are fulfilled.\\
{\rm (i)} $f$ is Hadamard differentiable at $\hat{x}$.\\
{\rm (ii)} For all $i\in \{1,...,p\}$, $g_i$ is Hadamard differentiable at $\hat{x}$ when $g_i(\hat{x})=0.$\\
{\rm (iii)} For all $i\in \{1,...,p\}$, $g_i$ is lower semicontinuous at $\hat{x}$ and G\^ateaux differentiable at $\hat{x}$ when $g_i(\hat{x})>0.$\\
{\rm (iv)} For all $j\in \{1,...,q\}$, $h_j$ is continuous on a neighborhood at $\hat{x}$ and Hadamard differentiable  at $\hat{x}$.\\
Then, there exist $\lambda_0,..., \lambda_p \in \R_+$ and $\mu_1,...,\mu_q \in \R$ which satisfy the following conditions.\\
\rm{(a)} $(\lambda_0,...,\, \lambda_p,\,\mu_1,...,\,\mu_q) \ne (0,...,0).$\\
\rm{(b)} For all $i \in \{1,..., p\},\; \lambda_ig_i(\hat{x})=0.$\\
\rm{(c)} $\lambda_0D_Hf(\hat{x})+\sum_{i=1}^{p} \lambda_iD_Gg_i(\hat{x})+\sum_{j=1}^{q} \mu_jD_Hh_j(\hat{x})=0.$\\
Futhermore, if we assume that the following assertion hold\\
{\rm (v)} $D_Hh_1(\hat{x}),...,\,D_Hh_q(\hat{x})$ are linearly independent\\
we have \\
{\rm (d)} $(\lambda_0,..., \lambda_p) \ne (0,...,0).$\\
Moreover, under (v) and the following assertion \\
{\rm (vi)} there exists $w\in \cap_{i=1}^{q} KerD_Hhj(\hat{x})$ such that for all $i\in \{1,...,p\}, \\D_Gg_i(\hat{x})w>0$ when $g_i(\hat{x})=0$\\
we can take \\
{\rm (e)} $\lambda_0=1$.\\
\end{theorem}
\section{Recall and Notations}
  
We set $\N$ the set of positive integer and $\N^*=\N\setminus\{0\}$. $\R$ denotes the set of real numbers and $\R_+$ the set of non negative real numbers.\\
Let $E$, $F$ and $G$ be three normed vector spaces, let $\Omega$ be a nonempty open subset of $E$, let $f: \Omega \rightarrow F$ be a mapping, let $x \in \Omega$, let $y\in E$ and $r \in ]0, +\infty[$. The closed ball centered at $y$ with a radius equal to $r$ is denoted by $ \overline{B}(y,r)$. \\
Let $A \subset E$ and $B\subset F$, $C^0(A,B)$ denotes the continuous mappings from $A$ into $B$. ${\rm bd} A$ denotes the topological boundary of $A$. \\
We denote by ${\mathcal L}(E,F)$ the space of the bounded linear mappings from $E$ into $F$.\\
Let $l\in {\mathcal L}(E,F)$, we note $Iml=l(E)$. Let $l_1\in {\mathcal L}(E,F)$  and $l_2 \in {\mathcal L}(E,G)$, we note by $(l_1,l_2)$ the mapping in ${\mathcal L}(E,F\times G)$  defined by for all $x\in E$, $(l_1,l_2)x=(l_1x,l_2x)$.    
$f$ is called G\^ateaux differentiable at $x$ when there exists $D_Gf(x)\in {\mathcal L}(X,Y)$ such that for all $h\in E$, $\lim_{t \downarrow 0} \frac{f(x+th)-f(x)}{t}=D_Gf(x)h$.\\
We say that $f$ is Hadamard differentiable at $x$ when there exists $D_Hf(x)\in {\mathcal L}(X,Y)$ such that for all $h\in E$, for all sequence $(h_n)_{n \in \N}$ converging to $h$ and for all sequence $(t_n)_{n \in \N}$ of positive numbers converging to $0$ we have\\ $\lim_{n \rightarrow +\infty} \frac{f(x+t_nh_n)-f(x)}{t_n}=D_Hf(x)h$ which is equivalent to (see \cite{F} p. 265) for each $K$ compact in $E$, $\lim_{t \downarrow 0} \sup_{h \in K}\frac{f(x+th)-f(x)}{t}=D_Hf(x)h$.\\
When $f$ is Hadamard differentiable at $x$, $f$ is also G\^ateaux differentiable at $x$ and $D_Hf(x)=D_Gf(x)$. But the converse is false when the dimension of E is greater than 2. \\
More information on these notions can be found in \cite{F}. \\ 
If $n\in \N^*$, we note $ \langle\cdot,\cdot \rangle$ the canonical scalar product on $\R^n$, $(e_{n,i})_{1\le i\le n}$ the canonical basis of $\R^n$ and $\|\cdot\|_{\infty}$ the maximum norm on $\R^n$. Moreover, we note by $\overline{B}_{\|\cdot\|_\infty}(y,r)$ the closed ball centered at $y$ with a radius equal to $r$ in $\R^n$ with the maximum norm. \\
We recall a supporting hyperplane theorem.
\begin{theorem}\label{th31} Let $n\in \N^*$. Let $C$ be a nonempty convex subset of $\R^n$ and $z \in {\rm bd}C$. Then there exist $v\in \R^n\setminus\{0\}$ and $\gamma \in \R$ such that $\langle v,z \rangle = \gamma$ and for all $x\in C$, $\langle v,x \rangle \le \gamma$.
\end{theorem}
\noindent
This theorem is a corralary of Hahn-Banach theorem. We can find a proof in \cite{JV} p. 37. Note that if $z=0$ we have $\gamma=0$.\\
\noindent
Next, we recall the Schauder fixed-point theorem.
\begin{theorem}\label{th32}(Schauder fixed-point theorem) Let $E$ be a normed vector space, let $C$ be a nonempty convex and compact subset of $E$ and let $f: C \rightarrow C$ be a continuous mapping, then $f$ admit a fixed point i.e. there exists $x\in C$ such that $f(x)=x$.\\
\end{theorem}
\noindent
We can find a proof of the Schauder fixed-point theorem in \cite{DG} p. 119.
\section{Proof of Theorem \ref{th21}}
We set $S:=\{i\in \{1,...,m\}:\; f_i(\hat{x})=0\}$. If $S=\emptyset$ we have for all $i\in \{1,..., m\},\; f_i(\hat{x})>0$ using the lower semicontinuous of $f_i$, there exists an open neighborhood $\Omega_1$ of $\hat{x}$ in $\Omega$ such that for all $i\in \{1,..., m\},$ for all $x \in \Omega_1$, $f_i(x)>0$. Since $\hat{x}$ is a solution of $({\mathcal I})$, we have $\hat{x}$ maximize  $f_0$ on $\Omega_1$. Therefore by using (i), we have $D_Gf_0(\hat{x})=0$. By taking $\lambda_0=1$ and for all $i\in \{1,..., m\}$, $\lambda_i=0$, we proved (a), (b), (c) and (d).
We assume that $S\ne \emptyset$ in the rest of the proof.
\subsection{To delete all inactive inequality constraints} 
By doing a change of index, we can assume that $S=\{1,..., s\}$ where $1\le s \le m$.
Since for all $i\in \{s+1,..., m\}$, we have $f_i(\hat{x})>0$, using (ii) there exists an open neighborhood $U$ of $\hat{x}$ in $\Omega$ such that for all $i\in \{s+1,..., m\}$, for all $x\in U$ , we have $f_i(x)>0$. For consequently, we have $\hat{x}$ is a solution of the following problem 
\[
({\mathcal NI})
\left\{
\begin{array}{cl}
{\rm Maximize} & f_0(x)  \\
{\rm subject \;  to} & x \in U \\
\null & \forall i\in\{1,..., s\},\, f_i(x) \ge 0.
\end{array}\right.
\] 
\subsection{ Proof of (a), (b), (c)} 
\vskip2mm
\noindent
We consider the mapping $F: U \rightarrow \R^{s+1}$ defined by $\forall x\in U$, $F(x)=(f_0(x),..., f_s(x))$. Since for all $i\in\{0,...,s\}$, $f_i$ is G\^ateaux differentiable at $\hat{x}$, we have $F$ is G\^ateaux differentiable at $\hat{x}$ and \\ $D_GF(\hat{x})=(D_Gf_0(\hat{x}),...,D_Gf_s(\hat{x})).$\\
We set $C:=ImD_GF(\hat{x})+\R_-^{s+1}$. We note that $C$ is a convex set of $\R^{s+1}$. 
Moreover, $C$ is not a neighborhood of 0. To prove this, we proceed by  contradiction, by assuming that $C$ is a neighborhood of 0. Therefore, there exists $r>0$ such that $\overline{B}_{\|\cdot\|_\infty}(0,r) \subset C$. Since $b=(r,..., r) \in \overline{B}_{\|\cdot\|_\infty}(0,r)$ we have $b\in C$ then there exists $u \in E$ and $z=(z_0,..., z_s)\in \R_-^{s+1}$ such that $D_GF(\hat{x})u + z=b$. For consequently, we have
\begin{equation}\label{eq1}
\left.
\begin{array}{r}
\forall i\in\{0,..., s\}, D_Gf_i(\hat{x}).u=r-z_i \ge r.
\end{array}
\right.
\end{equation}
By using (\ref{eq1}), we remark that $u\ne 0$.\\
Since $F$ is G\^ateaux differentiable at $\hat{x}$, we have 
\begin{equation}\label{eq2}
\left.
\begin{array}{r}
\exists \delta>0\; \forall t\in ]0,\delta]\;(\hat{x}+tu \in U) \; \|F(\hat{x}+tu)-F(\hat{x}) -tD_GF(\hat{x})u\|_\infty<rt.
\end{array}
\right.
\end{equation}
Then, using (\ref{eq2}) with $t=\delta$, we have $\|F(\hat{x}+\delta u)-F(\hat{x}) -\delta D_GF(\hat{x})u\|_\infty<r\delta$ which implies that $\forall i\in \{0,...,s\},f_i(\hat{x}+\delta u)-f_i(\hat{x}) -\delta D_Gf_i(\hat{x})u>-r\delta$. For consequently, by using (\ref{eq1}) we have for all  $i\in \{0,...,s\},f_i(\hat{x}+\delta u) -f_i(\hat{x})> \delta D_Gf_i(\hat{x})u-\delta r \ge 0$. Therefore, we have $f_0(\hat{x}+\delta u) > f_0(\hat{x})$ and for all  $i\in \{0,...,s\},\;f_i(\hat{x}+\delta u)>0$ which implies that $\hat{x}$ is not a solution of $({\mathcal NI})$. This is a contradiction. Since $0 \in C$ and $C$ is not a neighborhood of $0$, we have $0\in {\rm bd}C$.\\  
Since $C$ is a convex of $\R^{s+1}$ and $0\in {\rm bd}C$, by using Theorem \ref{th31} there exists $v=(\lambda_0,...,\, \lambda_s) \in \R^{s+1}\setminus\{0\}$ such that for all $x\in C$, $\langle v,x \rangle \le 0$. For consequently, we have 
\begin{equation}\label{eq3}
\left.
\begin{array}{r}
\forall u\in E,\; \forall z=(z_0,...,z_s)\in \R_-^{s+1},\; \sum_{i=0}^{s} \lambda_i(D_Gf_i(\hat{x})u+z_{i}) \le 0.
\end{array}
\right.
\end{equation}
We set for all $i \in \{s+1,...,\,m\}\; \lambda_{i}=0$.\\ Since $(\lambda_0,...,\, \lambda_s) \ne 0$, we have $(\lambda_0,...,\, \lambda_m) \ne 0$.
Let $i\in \{0,...,s\}$, by using (\ref{eq3}) with $u=0$ and $z=-e_{s+1,{i+1}}$, we have $-\lambda_i \le 0$ which implies that $\lambda_i \ge 0$. \\
We have also for all $i \in \{1,..., m\},\; \lambda_if_i(\hat{x})=0.$\\ 
By using (\ref{eq3}), with $z=0$, we have $\forall u\in E,\;   \sum_{i=0}^{s} \lambda_iD_Gf_i(\hat{x})u \le 0$ which implies that 
\begin{equation}\label{eq4}
\left.
\begin{array}{r}
\sum_{i=0}^{s} \lambda_iD_Gf_i(\hat{x})=0,
\end{array}
\right.
\end{equation}
therefore $\sum_{i=0}^{m} \lambda_iD_Gf_i(\hat{x})=0$.\\
We proved (a), (b) and (c).
\subsection{Proof of (d)}In addition, if we assume (iv), we have $\lambda_0 \ne 0$. We proceed by contradiction by assuming that $\lambda_0=0$. Since (iii) and $(\lambda_1,...,\lambda_s)\in\R_+^{s}\setminus\{0\}$, we have $\sum_{i=1}^{s}\lambda_iD_Gf_i(\hat{x})w>0$. By using (\ref{eq4}), we have $\sum_{i=1}^{s} \lambda_iD_Gf_i(\hat{x})w=0$. This a contradiction. Since, $\lambda_0 \ne 0$, by taking for all $i\in\{0,...,m\},\; \lambda_i'=\frac{\lambda_i}{\lambda_0}$, we proved (d).
\section{Proof of Theorem \ref{th22}}
We set $S:=\{i\in \{1,...,p\}:\; g_i(\hat{x})=0\}$. Without loss of generality, we can assume that $S \ne \emptyset$. If $S = \emptyset$ we can delete all inequality constraints. Indeed, we have for all $i\in\{1,...,\, p\}$, $g_i(\hat{x})>0$, by using (iii), there exists a neighborhood $\Omega_1$ of $\hat{x}$ in $\Omega$ such that for all $i\in\{1,...,\, p\}$, for all $x\in\Omega_1$, $g_i(x)>0$. For consequently, $\hat{x}$ is a solution of the following problem 
\[({\mathcal SP})
\left\{
\begin{array}{cl}
{\rm Maximize} & f(x)  \\
{\rm subject \;  to} & x \in \Omega_1 \\
\null & \forall j\in\{1,..., q\},\, h_j(x) = 0.\\
\end{array}\right.
\]
\subsection{To delete all inactive inequality constraints}
In the rest of the proof, we assume that $S \ne \emptyset$. By doing a change of index, we can assume that $S=\{1,..., s\}$ where $1\le s \le m$. Since for all $i\in \{s+1,..., m\}$, we have $g_i(\hat{x})>0$, using (iii) and (iv) there exists an open neighborhood $U$ of $\hat{x}$ in $\Omega$ such that for all $i\in \{s+1,..., m\}$, for all $x\in U$ , we have $g_i(x)>0$ and for all $j\in\{1,...,\,q\}$, $h_j$ is continuous on $U$. For consequently, we have $\hat{x}$ is a solution of the following problem 
\[
({\mathcal NP})
\left\{
\begin{array}{cl}
{\rm Maximize} & f(x)  \\
{\rm subject \;  to} & x \in U \\
\null & \forall i\in\{1,..., s\},\, g_i(x) \ge 0\\
\null & \forall j\in\{1,..., q\},\, h_j(x) = 0.\\
\end{array}\right.
\] 
\subsection{Proof of (a), (b), (c)}
We consider the mappings $G: U \rightarrow \R^{s+1}$ and $H: U \rightarrow \R^{q}$  defined by $\forall x\in U$, $G(x)=(f(x),g_1(x),..., g_s(x))$ and \\$H(x)=(h_1(x),...,h_q(x))$. Since (i), (ii) and (iv) we have $G$ and $H$ are Hadamard differentiable at $\hat{x}.$ Moreover $D_HG(\hat{x})=(D_Hf(\hat{x}), D_Hg_1(\hat{x})...,D_Hg_s(\hat{x}))$ and \\$D_HH(\hat{x})=(D_Hh_1(\hat{x}),...,D_Hh_q(\hat{x}))$.\\
We set $C:=Im(D_HG(\hat{x}),D_HH(\hat{x}))+\R_-^{s+1}\times \{0\}$. $C$ is a convex set of $\R^{s+q+1}$.\\ 
$C$ is not a neighborhood of 0. To prove this, we proceed by contradiction, by assuming that $C$ is a neighborhood of 0. Therefore, there exists $r>0$ such that $\overline{B}_{\|\cdot\|_\infty}(0,r) \subset C$.\\
We set $b=(r,...,r) \in \R^{s+1}$. 
Since, for all $j\in\{1,..., q\},$ $(b,re_{q,j})\in C$ and $(b,-re_{q,j})\in C$, there exists $u_j\in E$ and $z_j=(z_{0,j},...,\,z_{s,j}) \in \R_-^{s+1}$ such that
\begin{equation}\label{eq21}
\left.
\begin{array}{r}
D_HG(\hat{x})u_j +z_j = b \text{ and } D_HH(\hat{x})u_j=re_{q,j}
\end{array}
\right.
\end{equation} 
and 
there exists $\tilde{u}_j\in E$ and $\tilde{z}_j=(\tilde{z}_{0,j},...,\,\tilde{z}_{s,j}) \in \R_-^{s+1}$ such that
\begin{equation}\label{eq22}
\left.
\begin{array}{r}
D_HG(\hat{x})\tilde{u}_j +\tilde{z}_j = b \text{ and } D_HH(\hat{x})\tilde{u}_j=-re_{q,j}.
\end{array}
\right.
\end{equation} 
We set $K:=\{\sum_{j=1}^{q} a_ju_j + \sum_{j=1}^{q} \tilde{a}_j\tilde{u}_j\, :\, \forall j\in\{1,...,\,q\},\, a_j \ge 0, \tilde{a}_j \ge 0\text{ and }\\ \sum_{j=1}^{q} a_j + \sum_{j=1}^{q} \tilde{a}_j=1\}$. By using (\ref{eq21}), we have $0\notin K$. \\
We remark that $K$ is a convex and compact set of $E$. Since $G$ and $H$ are Hadamard differentiable at $\hat{x}$ and $K$ is a compact set of $E$, we have 
\begin{equation}\label{eq23}
\left.
\begin{array}{r}
\exists \delta_1>0\; \forall t\in ]0,\delta_1],\; \forall k\in K,\; \|G(\hat{x}+tk)-G(\hat{x}) -tD_HG(\hat{x})k\|_\infty<rt
\end{array}
\right.
\end{equation}
\begin{equation}\label{eq24}
\left.
\begin{array}{r}
\exists \delta_2>0\; \forall t\in ]0,\delta_2],\; \forall k\in K,\; \|H(\hat{x}+tk)-H(\hat{x}) -tD_HH(\hat{x})k\|_\infty<\frac{r}{q}t.
\end{array}
\right.
\end{equation}
Since $U$ is a neighborhood of $\hat{x}$, we have there exists $r_0 >0$ such that $ \overline{B}(\hat{x},r_0) \subset U$. we set $\alpha:=\min\{\delta_1,\delta_2,\frac{r_0}{\sup_{k \in K} \|k\|}\}. $
Therefore, by using (\ref{eq23}) and (\ref{eq24}) with $t=\alpha$, we have \\
\begin{equation}\label{eq25}
\left.
\begin{array}{r}
\forall k\in K,\, |f(\hat{x}+\alpha k)-f(\hat{x}) -\alpha D_Hf(\hat{x})k|<r\alpha
\end{array}
\right.
\end{equation}
\begin{equation}\label{eq26}
\left.
\begin{array}{r}
\forall k\in K,\,\forall i\in\{1,...,\,s\},\, |g_i(\hat{x}+\alpha k)-\alpha D_Hg_i(\hat{x})k|<r\alpha
\end{array}
\right.
\end{equation}
\begin{equation}\label{eq27}
\left.
\begin{array}{r}
\forall k\in K,\, \forall j\in\{1,...,\,q\},\, |h_j(\hat{x}+\alpha k)-\alpha D_Hh_j(\hat{x})k|<\frac{r}{q}\alpha.
\end{array}
\right.
\end{equation}
We set for all $j\in \{1,...,\, q\},\,$ for all $k\in K,\, w_j(k)=\frac{1}{\alpha}h_j(\hat{x}+\alpha k)-D_Hh_j(\hat{x})k.$
We note that for all $j\in \{1,...,\, q\},\,$ $w_j \in C^0(K,\R)$ because $h_j\in C^0(K,\R)$  and for all $j\in \{1,...,\, q\},\,$ for all $k\in K,\,|w_j(k)| <\frac{r}{q}$ (by (\ref{eq27})).\\
We consider the mapping $\Phi: K \rightarrow E$ defined by \\$\Phi(k):=\sum_{j=1}^{q} \phi_j(k)u_j+\sum_{j=1}^{q} \tilde{\phi}_j(k)\tilde{u}_j$ where for all $j\in \{1,...,q\},$ for all $k\in K$, $\phi_j(k)=\frac{1}{2q}-\frac{1}{2r}w_j(k)$ and $\tilde{\phi}_j(k)=\frac{1}{2q}+\frac{1}{2r}w_j(k)$. For all $k\in K$, we have for all $j\in \{1,...,q\},$ $\phi_j(k) \ge 0$, $\tilde{\phi}_j(k) \ge 0$ and $\sum_{j=1}^{q} \phi_j(k)+\sum_{j=1}^{q} \tilde{\phi}_j(k)=1$ therefore we have $\Phi(K) \subset K$.   
Since for all $j\in\{1,..., q\}$, $\phi_j$ and $\tilde{\phi}_j$ belong to $C^0(K,\R)$, we have $\Phi \in C^0(K,K)$. Therefore, by using Theorem \ref{th32}, we obtain there exists $\hat{k} \in K$ such that $\Phi(\hat{k})=\hat{k}$.
Since \[
\begin{array}{ll}
D_HH(\hat{x})\hat{k}&=D_HH(\hat{x})\Phi(\hat{k})\\
\null&=\sum_{j=1}^{q} (\frac{1}{2q}-\frac{1}{2r}w_j(\hat{k}))D_HH(\hat{x})u_j+\sum_{j=1}^{q} (\frac{1}{2q}+\frac{1}{2r}w_j(\hat{k}))D_HH(\hat{x})\tilde{u}_j \\
\null&=-\sum_{j=1}^{q}w_j(\hat{k})e_{q,j}. 
\end{array}
\]
Therefore, for all $j\in \{1,...,\,q\}, D_Hh_j(\hat{x})\hat{k}=-w_j(\hat{k})=-(\frac{1}{\alpha}h_j(\hat{x}+\alpha \hat{k})-D_Hh_j(\hat{x})\hat{k})$ which implies that $h_j(\hat{x}+\alpha \hat{k})=0.$\\
Since $\hat{k}\in K$ there exists $(a_j,\tilde{a}_j)_{1\le j\le q}\in \R_+^{2q}$ with $\sum_{j=1}^{q} a_j + \sum_{j=1}^{q} \tilde{a}_j=1$ such that $\hat{k}=\sum_{j=1}^{q} a_ju_j + \sum_{j=1}^{q} \tilde{a}_j\tilde{u}_j$. For consequently, we have \[
\begin{array}{ll}
D_Hf(\hat{x})\hat{k}&=\sum_{j=1}^{q} a_jD_Hf(\hat{x})u_j + \sum_{j=1}^{q} \tilde{a}_jD_Hf(\hat{x})\tilde{u}_j\\
\null&=\sum_{j=1}^{q} a_j(r-z_{0,j}) + \sum_{j=1}^{q} \tilde{a}_j(r-\tilde{z}_{0,j}) \text{ (from (\ref{eq21}) and (\ref{eq22}))}
\end{array}
\] which implies that
 \begin{equation}\label{eq28} 
\left.
\begin{array}{r}
D_Hf(\hat{x})\hat{k} \ge r.
\end{array}
\right.
\end{equation}
By the same reasoning, we have also 
 \begin{equation}\label{eq29} 
\left.
\begin{array}{r}
\forall i\in\{1,...,s\}, D_Hg_i(\hat{x})\hat{k} \ge r.
\end{array}
\right.
\end{equation}
By using (\ref{eq28}), (\ref{eq29}), (\ref{eq25})  and (\ref{eq26}) with $t=\hat{k}$, we have $f(\hat{x}+\alpha\hat{k})>f(\hat{x})$ and for all $i\in\{1,...,\,s\},$ we have $g_i(\hat{x}+ \alpha\hat{k})>0$.
Since $\hat{x}+\alpha\hat{k} \in U$ verify for all $i\in\{1,...,s\},$  $g_i(\hat{x}+ \alpha\hat{k})>0$, for all $j\in\{1,...,q\},$  $h_j(\hat{x}+ \alpha\hat{k})= 0$ and $f(\hat{x}+\alpha\hat{k})>f(\hat{x})$, we have $\hat{x}$ is not a solution of $(\mathcal{NP})$. This is a contradiction. Since $0 \in C$ and $C$ is not a neighborhood of $0$, we have $0\in {\rm bd}C$.\\ 
Since $C$ is a convex of $\R^{1+s+q}$ and $0\in {\rm bd}C$, by using Theorem \ref{th31} there exists $v=(\lambda_0,...,\, \lambda_s,\mu_1,...,\mu_q) \in \R^{1+s+q}\setminus\{0\}$ such that for all $x\in C$, $\langle v,x \rangle \le 0$. Therefore, we have 
\begin{equation}\label{eq30}
\left.
\begin{array}{l}
\forall u\in E,\; \forall z=(z_0,...,z_s)\in \R_-^{1+s}\\
\lambda_0(D_Hf(\hat{x})u+z_{0}) + \sum_{i=1}^{s} \lambda_i(D_Hg_i(\hat{x})u+z_{i})+\sum_{j=1}^{q} \mu_jD_Hh_j(\hat{x})u\le 0
\end{array}
\right\}
\end{equation}
We set for all $i \in \{s+1,...,\,p\}\; \lambda_{i}=0$. Since $(\lambda_0,...,\, \lambda_s,\mu_1,...,\mu_q)\ne 0$, we have $(\lambda_0,...,\, \lambda_p,\mu_1,...,\mu_q)\ne 0.$ \\ 
Let $i\in \{0,...,s\}$, by using (\ref{eq30}) with $u=0$ and $z=-e_{s+1,i+1}$, we have $-\lambda_i \le 0$ which implies that $\lambda_i \ge 0.$ 
We have also for all $i \in \{1,..., p\},\; \lambda_ig_i(\hat{x})=0.$
By using (\ref{eq30}), with $z=0$, we have \\ $\forall u\in E,\;   \lambda_0D_Hf(\hat{x})u + \sum_{i=1}^{s} \lambda_iD_Hg_i(\hat{x})u+\sum_{j=1}^{q} \mu_jD_Hh_j(\hat{x})u \le 0$ which implies that $\lambda_0D_Hf(\hat{x}) + \sum_{i=1}^{s} \lambda_iD_Hg_i(\hat{x})+\sum_{j=1}^{q} \mu_jD_Hh_j(\hat{x})=0.$ Since, for all $i\in\{1,...,s\}$, $D_Gg_i(\hat{x})=D_Hg_i(\hat{x})$, we have
\begin{equation}\label{eq31}
\left.
\begin{array}{r}
\lambda_0D_Hf(\hat{x})u + \sum_{i=1}^{s} \lambda_iD_Gg_i(\hat{x})u+\sum_{j=1}^{q} \mu_jD_Hh_j(\hat{x})u=0.
\end{array}
\right.
\end{equation}
Therefore $\lambda_0D_Hf(\hat{x})u + \sum_{i=1}^{p} \lambda_iD_Gg_i(\hat{x})u+\sum_{j=1}^{q} \mu_jD_Hh_j(\hat{x})u=0.$ We proved (a), (b) and (c).

\subsection{Proof of (d)}
We assume (i), (ii), (iii) and (iv). We proceed by contradiction by assuming that $(\lambda_0,..., \lambda_p) = (0,...,0)$. Therefore, according to (\ref{eq31}),\\ $\sum_{j=1}^{q} \mu_jD_Hh_j(\hat{x})u=0$. Since (iv), we have $(\mu_1,...,\mu_q)=0$. For consequently, we have $(\lambda_0,..., \lambda_p, \mu_1,...,\mu_q) = (0,...,0)$ this a contradiction with (a). We proved (d).  
\subsection{Proof of (e)}
We assume (i), (ii), (iii), (iv) and (v). Thanks to our previous proof we know that there exist  $\lambda_0,..., \lambda_p \in \R_+$ and $\mu_1,...,\mu_q \in \R$ which verify (a), (b), (c) and (d). We have $\lambda_0 \ne 0$, we proceed by contradiction by assuming that $\lambda_0=0$. 
Since (d) and (b), we have $(\lambda_1,...,\,\lambda_s)\ne 0$. Since (v) and $(\lambda_1,...,\,\lambda_s)\in \R_+^s\setminus\{0\}$, we have $\sum_{i=1}^{s} \lambda_iD_Gg_i(\hat{x})w >0$. By using (\ref{eq31}), we have $\sum_{i=1}^{s} \lambda_iD_Gg_i(\hat{x})w=0$. This a contradiction. Since, $\lambda_0 \ne 0$, by taking for all $i\in\{0,...,p\},\; \lambda_i'=\frac{\lambda_i}{\lambda_0}$, we proved (e).   
\vskip2mm
\noindent
{\textbf{ Acknowledgements}} The author thanks Jo\"el Blot for his encouragement.

\end{document}